\documentclass[preprint,11pt]{article}

\usepackage{epstopdf}
\usepackage[caption=false]{subfig}

\usepackage[numbers,sort&compress]{natbib}
\bibpunct[, ]{[}{]}{,}{n}{,}{,}
\renewcommand\bibfont{\fontsize{10}{12}\selectfont}
\makeatletter
\def\NAT@def@citea{\def\@citea{\NAT@separator}}
\makeatother

\usepackage{comment}
\usepackage{amsmath,amssymb,amsthm}
\usepackage{url}
\usepackage{color}
\usepackage{ulem}
\usepackage{graphicx}
\usepackage{bbm}
\usepackage{booktabs}
\usepackage{multirow}
\usepackage{natbib}
\usepackage{authblk}

\theoremstyle{plain}
\newtheorem{theorem}{Theorem}[section]
\newtheorem{lemma}[theorem]{Lemma}

\theoremstyle{definition}

\theoremstyle{remark}

\newcommand{\N}{\mathbb{N}}

\newcommand{\R}{\mathbb{R}}

\newcommand{\bX}{\boldsymbol{X}}

\newcommand{\bvt}{\boldsymbol{\vartheta}}

\begin{document}


\title{Confidence bands for exponential distribution functions under progressive type-II censoring}

\author{
Stefan Bedbur\thanks{Corresponding author: Stefan Bedbur, bedbur@isw.rwth-aachen.de} \qquad Fabian Mies}
\affil{RWTH Aachen University, Germany}

\date{December 22, 2020}

\maketitle

\begin{abstract}
Based on a progressively type-II censored sample from the exponential distribution with unknown location and scale parameter, confidence bands are proposed for the underlying distribution function by using confidence regions for the parameters and Kolmogorov-Smirnov type statistics. Simple explicit representations for the boundaries and for the coverage probabilities of the confidence bands are analytically derived, and the performance of the bands is compared in terms of band width and area by means of a data example. As a by-product, a novel confidence region for the location-scale parameter is obtained. Extensions of the results to related models for ordered data, such as sequential order statistics, as well as to other underlying location-scale families of distributions are discussed.

\textbf{Keywords:} confidence set; coverage probability; finite sample inference; location-scale family; order statistics; simultaneous confidence intervals
\end{abstract}

\section{Introduction}
In life testing, progressively type-II censored order statistics serve to model component lifetimes in experiments, where upon each failure a pre-fixed number of intact components are removed; see, e.g., \cite{BalAgg2000,BalCra2014}. Such withdrawals of operating components may be part of the experimental design for different reasons. As two examples, a researcher may wish to release capacities for another experiment or to save costs in situations, where the ongoing operation of the components is expensive. For more motivational aspects on the model, we refer to \cite{Bal2007,Cra2017}.   

A progressively type-II censored lifetime experiment is formally described as follows. Suppose that we put $n$ identical components on a lifetime test. At time of the $j$-th failure of some component, $1\leq j\leq m$, $R_j\in\N_0$ of the still operating components are (randomly) selected and removed from the experiment. Upon the $m$-th component failure, all remaining $R_m=n-m-\sum_{j=1}^{m-1}R_j$ components are removed. Thus, $m$ component failure times are recorded while $\sum_{j=1}^m R_j=n-m$ component failure times have been (progressively type-II) censored. 
Note that conventional type-II right censoring is included in this setup by choosing $R_1=\dots=R_{m-1}=0$ and $R_m=n-m$. In that case, choosing $m=n$ corresponds to a complete sample from the underlying lifetime distribution.

Inferential results for the underlying lifetime distribution of progressively type-II censored order statistics have widely been studied. We refer to \cite{BalCra2014} for a summary providing many references; see also \cite{Cra2017}. Among others, confidence regions have been proposed for parameters, quantiles, and (point-wise) reliability of the lifetime distribution in different distribution families including exponential, Weibull, and Pareto models. Recent works on inference under progressive type-II censoring are provided, e.g., by \cite{BalHayLiuKia2015,DoeCra2019,MalKha2020} for the exponential distribution, and by \cite{BaiShiLiuLiu2019,MonKun2019} for other parametric distributions.

The scope of the present work is to develop, study, and compare confidence bands for the underlying cumulative distribution function (cdf) of progressively type-II censored order statistics. Such confidence bands are random sets in the two-dimensional plane, which contain the entire graph of the true cdf with a pre-specified probability. The bands may be interpreted as simultaneous confidence intervals around estimated survival probabilities, effectively addressing the multiple-testing problem. This is of particular relevance for applications in reliability, where a statistical statement for the failure probabilities of some unit jointly at several time points is required. 

Asymptotic confidence bands may be constructed based on a nonparametric estimator for the survival function (or cdf); see, e.g., \cite{Nai1984,HolMckYan1997}. However, in life testing, the acquisition of samples is typically costly, and small sample sizes may render asymptotic methods inapplicable while parametric models still allow for exact inference (provided they are properly specified). In this paper, we focus on a single sample from the two-parameter exponential distribution although the methodologies presented here also apply to the multi-sample case and other location-scale families. Explicit formulas are presented for various confidence bands satisfying a desired exact confidence level. Two of these confidence bands ($B_1$ and $B_4$) have been studied by \cite{SriKanWha1975} for the case of a complete sample of independent and identically distributed (iid) exponential random variables; here, the results are extended to a progressively type-II censored sample. Moreover, we show that the confidence band $B_4$ may be trimmed without reducing its coverage probability. Beyond that, we propose two more confidence bands ($B_2$ and $B_3$) and discuss their properties. As a by-product, our methodology also suggests a new explicit confidence region ($C''_4$) for the location-scale parameter of the exponential distribution. 

The detailed outline of the paper is as follows. 
In Section \ref{s:PCII}, the model of progressively type-II censored order statistics is introduced, and some properties are stated.  Based on a sample of progressively type-II censored order statistics from a two-parameter exponential distribution,  confidence bands for the underlying cdf are then proposed in Section \ref{s:bands}. Here, confidence bands are constructed via confidence regions for the parameters in Section \ref{ss:conf} and via Kolomogorov-Smirnov type statistics in Section \ref{ss:kst}. Explicit formulas for the boundaries and the coverage probabilities of the confidence bands are presented, and it is shown how to obtain the quantiles of the relevant distributions by simulation. In Section \ref{s:comp}, the bands are illustrated and compared with respect to band width and area for a classical data set. Extensions of the results are discussed in Section \ref{s:general}, where we focus on related models for ordered data in Section \ref{ss:sos} and consider other underlying location-scale families in Section \ref{ss:otherlsf}. Section \ref{s:conclusion} finally gives the conclusion highlighting the main findings.


\section{Progressively type-II censored order statistics}\label{s:PCII}

For $m,n\in\N$ with $1<m\leq n$, we assume to have a sample $X_{1:m:n}\leq\dots\leq X_{m:m:n}$ of progressively type-II censored order statistics based on some absolutely continuous cdf $F$ with corresponding density function $f$ and the censoring scheme $(R_1,\dots,R_m)\in\N_0^m$ with $\sum_{j=1}^m R_j=n-m$. The joint density function of $X_{1:m:n},\dots, X_{m:m:n}$ is then given by 
\begin{equation*}
	f(x_1,\dots,x_m)\,=\,\prod_{j=1}^m \gamma_jf(x_j)[1-F(x_j)]^{R_j}
\end{equation*}
for $x_1\leq\dots\leq x_m$, 
where 
\begin{equation}\label{eq:gammas}
	\gamma_j\,=\,\sum_{i=j}^m (R_i+1)\,,\quad 1\leq j\leq m\,,
\end{equation} denote (known) positive parameters; see, e.g., \cite{BalCra2014}, p. 22. Note that $\gamma_1=n$. To shorten notation, let $\bX=(X_{1:m:n},\dots,X_{m:m:n})$ and $\boldsymbol{R}=(R_1,\dots,R_m)$.

Now, let us assume that the underlying cdf $F$ belongs to the location-scale family $\mathcal{F}=\{F_{\bvt}:\bvt\in\Theta\}$, $\Theta=\R\times(0,\infty)$, of exponential cdfs, i.e.,
\begin{equation}\label{eq:locscale}
	F_{\bvt}(x)\,=\,1-\exp\left\{-\,\frac{x-\mu}{\sigma}\right\}\,,\quad x>\mu\,,
\end{equation}
for $\bvt=(\mu,\sigma)\in\Theta$. Throughout this article, the location-scale parameter $\bvt$ is supposed to be unknown. It is well known that the maximum likelihood estimator (MLE) $\hat{\bvt}=(\hat{\mu},\hat{\sigma})$ of $\bvt$ based on $\bX$ is given by
\begin{eqnarray}
	\hat{\mu}\,&=&\ X_{1:m:n}\label{eq:mudach}\\[1ex]
	\text{and}\quad\hat{\sigma}\,&=&\,\frac{1}{m}
	\sum_{j=2}^m \gamma_j(X_{j:m:n}-X_{j-1:m:n})\,,\label{eq:sigmadach}
\end{eqnarray}
where $\hat{\mu}$ and $\hat{\sigma}$ are independent with $\hat{\mu}\thicksim F_{(\mu,\sigma/n)}\in\mathcal{F}$ and $\hat{\sigma}\thicksim\Gamma(m-1,\sigma/m)$. Here, $\Gamma(a,b)$ denotes the gamma distribution with shape parameter $a>0$, scale parameter $b>0$, and mean $ab$. Moreover, $\hat{\bvt}$ is a complete sufficient statistic for $\bvt$, and the uniformly minimum variance unbiased estimators of $\mu$ and $\sigma$ are given by $\tilde{\mu}=\hat{\mu}-\tilde{\sigma}/n$ and $\tilde{\sigma}=m\hat{\sigma}/(m-1)$; for more details, see \cite{BalCra2014}, Section 12.





\section{Confidence bands for the baseline cdf}\label{s:bands}

Based on progressively type-II censored order statistics $\bX$ with censoring scheme $\boldsymbol{R}$, we aim for constructing confidence bands for (the graph) of the underlying cdf $F_{\bvt}$ subject to a desired (exact) confidence level. For this, let $\bX$ be formally defined on the probability space $(\Omega,\mathcal{A},P_{\bvt})$. A confidence band is then introduced as a random mapping $B=B(\bX)$ with values in the power set of $\R^2$, which satisfies the following two properties for all $\bvt\in\Theta$:
\begin{itemize}
	\item[(i)] $\{\text{graph}\,F_{\bvt}\subseteq B\}\in\mathcal{A}$, where $\text{graph}\, F_{\bvt}=\{(t,F_{\bvt}(t)):t\in\R\}$ denotes the graph of $F_{\bvt}$, and
	\item[(ii)] the sets $\{y:(t,y)\in B\}$, $t\in\R$, form (possibly degenerated) intervals $P_{\bvt}$-almost-surely. 
\end{itemize}
A confidence band $B$ is then said to have (exact) confidence level $1-p\in(0,1)$ if for all $\bvt\in\Theta$
\begin{equation*}
	P_{\bvt}(\text{graph}\,F_{\bvt}\subseteq B)\stackrel{(=)}{\geq} 1-p\,.
\end{equation*}
For the construction of confidence bands for $F_{\bvt}$, we focus on two parametric methods proposed in the literature. The first approach is presented by \cite{Kan1968b,Kan1968a} and based on the availability of some confidence region for $\bvt$;  see also \cite{MieBed2017}. The second method is developed by \cite{KanSri1972} and makes use of so called Kolmogorov-Smirnov type statistics, i.e., parametric analogues of the nonparametric Kolmogorov-Smirnov statistic.

Finally, note that once having derived confidence bands for the underlying cdf, confidence bands for the corresponding reliability function and quantile function can be obtained by respective transformations. For example, if $B$ denotes a confidence band for $F_{\bvt}$ with (exact) confidence level $1-p$, then
\begin{equation*}
	\overline{B}\,=\,\{(x,1-y):(x,y)\in B\}
\end{equation*}
forms a confidence band for the reliability function $1-F_{\bvt}$ with (exact) confidence level $1-p$.

\subsection{Confidence Bands based on Confidence Regions}\label{ss:conf}

Whenever there is some confidence region for $\bvt$ available, one may proceed as follows to obtain a confidence band for $F_{\bvt}$; see \cite{Kan1968b,Kan1968a}. Let $p\in(0,1)$ and $C=C(\bX)$ be a confidence region for $\bvt$ with exact confidence level $1-p$. Moreover, we assume that $C$ is path-connected $P_{\bvt}$-almost-surely for every $\bvt\in\Theta$. Then,
\[B_C\,=\,\bigcup_{\tilde{\bvt}\in C}\text{graph}\,F_{\tilde{\bvt}}\]
forms a confidence band for $F_{\bvt}$ provided that it meets the measurability condition (i). If so, it is obvious from the definition that $B_C$ has at least confidence level $1-p$. However, to prevent large and thus less informative bands, an exact confidence level of $1-p$ may be desired. For this, a sufficient condition is that $C$ satisfies the implication
\begin{equation}
	\text{graph}\,F_{\tilde{\bvt}}\subseteq B_C\quad\Rightarrow\quad \tilde{\bvt}\in C \label{eq:exhaustive}
\end{equation} 
$P_{\bvt}$-almost-surely for every $\bvt\in\Theta$. If condition \eqref{eq:exhaustive} holds, we say that the set $C$ is exhaustive. In that case, the measurability condition (i) on $B_C$ is trivially met, and we have equal coverage probabilities 
\begin{equation*}
	P_{\bvt}(\text{graph}\,F_{\tilde{\bvt}}\subseteq B_C)\,=\,P_{\bvt}(\tilde{\bvt}\in C),\quad \tilde{\bvt}\in\Theta\,,
\end{equation*}
for every $\bvt\in\Theta$; in particular, $B_C$ is then unbiased if and only if $C$ is unbiased. 

Based on an iid sample from the general location-scale family, characterizations of exhaustive confidence regions are provided in \cite{MieBed2017} by making a case distinction on the supports of the underlying cdfs; the results therein are readily seen to remain true for a progressively type-II censored sample. In the present situation, we have left-bounded supports and thus conclude from the results in \cite{MieBed2017} that any compact (and path-connected) confidence region is exhaustive if and only if it is convex and comprehensive. Here, a set $A\subseteq\R^2$ is called comprehensive if  $(x_1,x_2),(y_1,y_2)\in A$ and $(z_1,z_2)\in\R^2$ with $x_i\leq z_i\leq y_i$, $i=1,2$, imply that $(z_1,z_2)\in A$.

Two different confidence regions for $\bvt$ have been constructed in \cite{Wu2010} by combining the independent pivotal quantities
\begin{eqnarray*}
	\frac{n(\hat{\mu}-\mu)}{\sigma}\thicksim F_{(0,1)}\,,\qquad \frac{2m\hat{\sigma}}{\sigma}\thicksim \chi^2(2m-2)
\end{eqnarray*}
and
\begin{eqnarray*}
	\frac{n(\hat{\mu}-\mu)}{m\hat{\sigma}/(m-1)}\thicksim \text{F}(2,2m-2)\,,\qquad
	\frac{2(n(\hat{\mu}-\mu)+m\hat{\sigma})}{\sigma}\thicksim \chi^2(2m)\,,
\end{eqnarray*}
respectively. Here, $\chi^2(k)$  and $\text{F}(k_1,k_2)$ denote the chi-square distribution and F-distribution with $k$ and $k_1,k_2$ degrees of freedom, and $F_{(0,1)}\in\mathcal{F}$ is the standard exponential distribution, see formula \eqref{eq:locscale}. Allocating the overall confidence level uniformly to the intervals and tails, confidence regions with exact confidence level $1-p\in(0,1)$ are thus
\begin{eqnarray}
	C_1\,=\,\Big\{(\mu,\sigma)\in\Theta\,:\,a_{q_1}(\sigma)\leq\mu\leq a_{q_2}(\sigma)\,,
	\, \sigma_{q_2}\leq\sigma\leq\sigma_{q_1}\Big\}\label{eq:Wu1}
\end{eqnarray}
with
\begin{eqnarray*}
	a_q(\sigma)\,&=&\,\hat{\mu}+\frac{\sigma\ln(q)}{n}\,,\quad\sigma>0\,,\\[1ex]
	\text{and}\qquad \sigma_q\,&=&\,\frac{2m\hat{\sigma}}{\chi_{q}^2(2m-2)}\,,\quad q\in(0,1)\,,
\end{eqnarray*}
as well as
\begin{eqnarray}
	C_2\,=\,\Big\{(\mu,\sigma)\in\Theta\,:\,\mu_{q_2}\leq\mu\leq\mu_{q_1}\,,\,b_{q_2}(\mu)\leq\sigma\leq b_{q_1}(\mu)\Big\}\label{eq:Wu2}
\end{eqnarray}
with
\begin{eqnarray*}
	b_q(\mu)\,&=&\,\frac{2(n(\hat{\mu}-\mu)+m\hat{\sigma})}{\chi_{q}^2(2m)}\,,\quad\mu\in\R\,,\\[1ex]
	\text{and}\qquad
	\mu_q\,&=&\,\hat{\mu}-\frac{m\hat{\sigma}\text{F}_{q}(2,2m-2)}{(m-1)n}\,,\quad q\in(0,1)\,.
\end{eqnarray*}
Here, $\chi^2_\beta(k)$ and $\text{F}_{\beta}(k_1,k_2)$ denote the $\beta$-quantile of $\chi^2(k)$ and $\text{F}(k_1,k_2)$, respectively, and we set $q_1=(1-\sqrt{1-p})/2$ and $q_2=1-q_1$, for brevity.

The confidence regions $C_1$ and $C_2$ have trapezoidal shape in the $(\mu,\sigma)$-plane; see Figures \ref{fig:C1} and \ref{fig:C2} for an illustration. $C_1$ is the area enclosed by the horizontal lines
$\sigma=\sigma_{q_1}$ and $\sigma=\sigma_{q_2}$ parallel to the $\mu$-axis and the diagonal lines
\[\sigma=\frac{n(\mu-\hat{\mu})}{\ln(q_i)}\,,\qquad i=1,2\,,\]
with negative slopes. Similarly, $C_2$ is bounded by the vertical lines $\mu=\mu_{q_1}$ and $\mu=\mu_{q_2}$ parallel to the $\sigma$-axis and by the diagonal lines
$\sigma=b_{q_1}(\mu)$, $\mu\in\R$, and $\sigma=b_{q_2}(\mu)$, $\mu\in\R$,
with negative slopes. Hence, $C_1$ and $C_2$ are both comprehensive and thus exhaustive; see \cite{MieBed2017}. The corresponding confidence bands
\begin{eqnarray*}
	B_i\,&\equiv&\, B_{C_i}\,=\,\bigcup_{\tilde{\bvt}\in C_i}\text{graph}\,F_{\tilde{\bvt}}\,,\qquad i=1,2\,,
\end{eqnarray*}
for $F_{\bvt}$ therefore have exact confidence level $1-p$ as it is the case for the underlying confidence regions. 

To have explicit representations for $B_1$ and $B_2$ at hand, we shall state the lower and upper boundaries of both confidence bands. For $B_1$, their derivation is in analogy to that in \cite{SriKanWha1975}, Section 2.2, for a complete exponential sample and therefore omitted; cf. also \cite{Hay2012}. In case of $B_2$, the proof can similarly be performed and is presented in the appendix.

\begin{theorem}\label{thm:Wu}
	Let the confidence regions $C_1$ and $C_2$ in formulas (\ref{eq:Wu1}) and (\ref{eq:Wu2}) have exact confidence level $1-p\in(0,1)$.
	Then,
	\begin{eqnarray*}
		B_i\,=\,\{(x,y)\in\R\times[0,1]:\,U_i(x)\leq y\leq O_i(x)\}
		,\quad i=1,2\,,
	\end{eqnarray*}
	with
	\begin{eqnarray*}
		U_1\,(x)&=&\begin{cases} F_{(a_{q_2}(\sigma_{q_2}),\sigma_{q_2})}(x)\,, & x\leq \hat{\mu}\\
			F_{(a_{q_2}(\sigma_{q_1}),\sigma_{q_1})}(x)\,, & x>\hat{\mu}\end{cases}\,,\\[1ex]
		O_1(x)\,&=&\begin{cases} F_{(a_{q_1}(\sigma_{q_1}),\sigma_{q_1})}(x)\,, & x\leq \hat{\mu}\\
			F_{(a_{q_1}(\sigma_{q_2}),\sigma_{q_2})}(x)\,, & x>\hat{\mu}\end{cases}\,,
	\end{eqnarray*}
	respectively
	\begin{eqnarray*}
		U_2(x)\,&=&\begin{cases} F_{(\mu_{q_1},b_{q_1}(\mu_{q_1}))}(x)\,, & x\leq \hat{\mu}+m\hat{\sigma}/n\\
			F_{(\mu_{q_2},b_{q_1}(\mu_{q_2}))}(x)\,, & x>\hat{\mu}+m\hat{\sigma}/n\end{cases}\,,\\[1ex]
		O_2(x)\,&=&\begin{cases} F_{(\mu_{q_2},b_{q_2}(\mu_{q_2}))}(x)\,, & x\leq \hat{\mu}+m\hat{\sigma}/n\\
			F_{(\mu_{q_1},b_{q_2}(\mu_{q_1}))}(x)\,, & x>\hat{\mu}+m\hat{\sigma}/n\end{cases}\,,
	\end{eqnarray*} 
	form confidence bands for $F_{\bvt}$ with exact confidence level $1-p$. 
\end{theorem}

In recent years, confidence regions with smallest area for distribution parameters of progressively type-II censored order statistics have been proposed; see, for instance, \cite{Fer2014} and \cite{AsgFerAbd2017} for the Pareto distribution and Rayleigh distribution. In \cite{BedKamLen2019}, a minimum area confidence region for an underlying location-scale parameter has been derived based on independent progressively type-II censored samples. For the one-sample exponential case, the finding yields that
\begin{eqnarray}
	C_3\,&=&\,\Bigg\{(\mu,\sigma)\in(-\infty,\hat{\mu}]\times(0,\infty):\notag\\[1ex]
	\quad&&(m+1)\ln\left(\frac{\hat{\sigma}}{\sigma}\right)
	\,-\,\frac{n(\hat{\mu}-\mu)+m\hat{\sigma}}{\sigma}\,\geq\, c_p\Bigg\}\label{eq:Cmin}
\end{eqnarray}
has smallest area among all confidence regions for $\bvt$ with confidence level $1-p\in(0,1)$, which are based on the pivotal quantity $((\hat{\mu}-\mu)/\sigma,\hat{\sigma}/\sigma)$; for the iid case, see also \cite{Zha2018}. Here, $c_p\equiv c_p(m)$ denotes the $p$-quantile of the distribution of the random variable
\[(m+1)\ln(Y)-mY-Z\,,\]
where $Y\thicksim\Gamma(m-1,1/m)$ and $Z\thicksim F_{(0,1)}\in\mathcal{F}$ are independent; see Section \ref{s:PCII}. In this context, a confidence region $C=C(\bX)$ is said to be based on some pivotal quantity $\boldsymbol{T}=\boldsymbol{T}(\bvt,\bX)$ if there exists a Borel set $A\subseteq\R^2$ with the property that $\bvt\in C(\bX)$ if and only if $\boldsymbol{T}(\bvt,\bX)\in A$. It is evident that the class of  all confidence regions based on $\boldsymbol{T}$ is invariant under measurable bijective transformations of $\boldsymbol{T}$. Hence, in the present situation, $C_3$ is seen to have smaller area than $C_1$ and $C_2$. Note that, in general, this relation does not necessarily transfer to the area of the corresponding confidence bands.

The algebraic structure of $C_3$ is seen to be the same as in the particular case of type-II right censoring; see formula (4) in \cite{LenBedKam2019}, Section 3. By inspecting the proof of Theorem 2 in \cite{LenBedKam2019}, we therefore find a more explicit representation for $C_3$, i.e.,
\begin{equation}\label{eq:C3g}
	C_3\,=\,\left\{(\mu,\sigma)\in\Theta:\hat{\mu}+g(\sigma)\leq\mu\leq\hat{\mu}\,,\,Z_{-1}\leq\sigma\leq Z_0\right\}\
\end{equation}
with mapping $g:(0,\infty)\rightarrow\R$ defined by
\begin{equation}\label{eq:g}
	g(z)\,=\,\frac{[c_p-(m+1)(\ln(\hat{\sigma})-\ln(z))]\,z+m\hat{\sigma}}{n}
\end{equation}
for $z>0$, and
\begin{equation*}
	Z_i\,=\,-\,\frac{m\hat{\sigma}}{m+1}\,\left[ W_i\left(-\,\frac{m}{m+1}\exp\left\{\frac{c_p}{m+1}\right\}\right)\right]^{-1}
\end{equation*}
for $i\in\{-1,0\}$. Here, $W_{-1}$ and $W_0$ denote the real-valued branches of the Lambert W-function, i.e., $W_{-1}$ $(W_0)$ is the inverse function of the mapping $U(x)=x\exp\{x\}$, $x\leq-1\;(x\geq-1)$.

The corresponding confidence band $B_3\equiv B_{C_3}$ is now as follows; for the derivation, see the appendix. 

\begin{theorem}\label{thm:minvol}
	Let the confidence region $C_3$ in formulas (\ref{eq:Cmin}) and (\ref{eq:C3g}) have exact confidence level $1-p\in(0,1)$.
	Then, $B_3$ has confidence level $1-p$ and is given by
	\begin{eqnarray*}
		B_3\,=\,\{(x,y)\in\R\times[0,1]:\,U_3(x)\leq y\leq O_3(x)\}
		\,,
	\end{eqnarray*}
	with
	\begin{equation*}
		U_3(x)\,=\,F_{(\hat{\mu},Z_0)}(x)\,,\quad  x\in\R\,,
	\end{equation*}
	and
	\begin{eqnarray*}
		O_3(x)\,=\,\begin{cases} F_{(\hat{\mu},Z_{-1})}(x)\,, & \sigma_x^*<Z_{-1}\\
			F_{(\hat{\mu}+g(\sigma_x^*),\sigma_x^*)}(x)\,, & Z_{-1}\leq \sigma_x^*\leq Z_0\\
			F_{(\hat{\mu},Z_0)}(x)\,, & \sigma_x^*>Z_0\end{cases}\,,
	\end{eqnarray*} 
	where
	\begin{equation*}
		\sigma_x^*\,=\,\frac{n(\hat{\mu}-x)+m\hat{\sigma}}{m+1}\,.
	\end{equation*}
\end{theorem}

The minimum area confidence region is depicted in Figure \ref{fig:C3}; see also \cite{LenBedKam2019} for the doubly type-II censored case. These figures indicate that $C_3$ is  not comprehensive and, hence, not exhaustive. That is, there exist parameters $\tilde{\bvt}\notin C_3$ satisfying $\text{graph}\, F_{\tilde{\bvt}}\subseteq B_3$, and an example of this is highlighted in Figure \ref{fig:C3}. The confidence band $B_3$ therefore has a confidence level greater than that of $C_3$. More precisely, the exact confidence level of $B_3$ is given by that of the comprehensive convex hull of $C_3$, which is defined as the smallest comprehensive and convex confidence region containing $C_3$. For an illustration, see again Figure \ref{fig:C3}; cf. \cite{MieBed2017}. The coverage probability of this superset can be computed analytically; the derivation can be found in the appendix.

\begin{lemma}\label{La:CPMin}
	Under the assumptions of Theorem \ref{thm:minvol}, the exact confidence level of $B_3$ is given by
	\begin{eqnarray*}
		\tau\,&=&\,1-p+\frac{e^{c_p}m^{m+1}}{2(m-2)!}\left(\frac{1}{y^2}-\frac{1}{z^2}\right)
		+\left(\frac{z}{m+1}\right)^{m-1}\\[1ex]
		&&\quad\times\left[G_{m-1}\left(\frac{(m+1)y}{z}\right)-G_{m-1}(m+1)\right]\,,
	\end{eqnarray*}
	where
	\begin{eqnarray*}
		y\,&=&\,-(m+1)W_0\left(-\,\frac{m}{m+1}\exp\left\{\frac{c_p}{m+1}\right\}\right)\,,\\[1ex]
		z\,&=&\,m\exp\left\{1+\frac{c_p}{m+1}\right\}\,,
	\end{eqnarray*}
	and $G_k$ denotes the cdf of $\Gamma(k,1)$. In particular, it does only depend on $p$ and $m$.
\end{lemma}

Table \ref{table:B3CPTab} shows the exact confidence level $\tau$ of $B_3$ for different values of $m$, where the exact confidence level $1-p$ of the underlying confidence region $C_3$ is chosen as 90\% and 95\%, respectively. Here, the constant $c_p$ in formula (\ref{eq:Cmin}) is each numerically obtained by a Monte Carlo simulation of size $10^7$.

\begin{table}[h!]
	\centering		
	\begin{tabular}{cc|rrrrrrrr}
		\toprule
		&\multirow{2}{*}{$1-p$}  &\multicolumn{8}{|c}{$m$}\\[1ex]
		&& 2 & 3 & 4 & 5 & 10  & 25 & 50  & 100
		\\ \midrule
		&$90\%$  & 91.1 & 91.7 & 92.0 & 92.2 & 92.5 & 92.6 & 92.6 & 92.5\\[1ex]
		&$95\%$ &  95.6 & 95.9 & 96.1 & 96.2 & 96.5 & 96.5 & 96.5 & 96.4 \\ \midrule
	\end{tabular}		
	\caption{Exact confidence level $\tau$ of $B_3$ (in \%) for exact confidence level $1-p$ of $C_3$ and different values of $m$.}	
	\label{table:B3CPTab}
\end{table}

Likewise, we may invert the formula in Lemma \ref{La:CPMin} numerically and choose $p=p(\tau)$ in such a way that $B_3$ meets some desired exact confidence level $\tau\in(0,1)$. The resulting critical value $c_{p(\tau)}$ is presented in Table \ref{table:B3Quantiles} for $\tau\in\{90\%,95\%\}$ and values of $m$ as in Table \ref{table:B3CPTab}. Here, we use a Monte Carlo simulation to approximate the quantile function $p\mapsto c_p$.

\begin{table*}[h!]
	\centering	
	\footnotesize	
	\begin{tabular}{cc|rrrrrrrr}
		\toprule 
		& \multirow{2}{*}{$\tau$} & \multicolumn{8}{c}{$m$}\\[1ex]
		& &  2 & 3 & 4 & 5 & 10 & 25 & 50 & 100 \\ 
		\midrule
		$1-p(\tau)$ & \multirow{2}{*}{90\%} & 88.8 & 88.0 & 87.6 & 87.4 & 86.9 & 86.8 & 86.9 & 87.0 \\ 
		$c_{p(\tau)}$ & & -9.784 & -8.372 & -8.542 & -9.116 & -13.385 & -28.025 & -52.924 & -102.878 \\ \midrule
		$1-p(\tau)$ & \multirow{2}{*}{95\%} & 94.4 & 93.9 & 93.6 & 93.4 & 93.1 & 93.0 & 93.1 & 93.1 \\ 
		$c_{p(\tau)}$ & & -11.906 & -9.807 & -9.737 & -10.191 & -14.272 & -28.806 & -53.684 & -103.614 \\ 
		\midrule	
	\end{tabular}
	\caption{Exact confidence level $1-p(\tau)$ (in \%) and critical value $c_{p(\tau)}$ to choose for $C_3$ such that $B_3$ has exact confidence level $\tau$ (obtained by $10^7$ Monte Carlo simulations, each).}	
	\label{table:B3Quantiles}
\end{table*}

\subsection{Kolomogorov-Smirnov Type Bands}\label{ss:kst}

We now focus on another method for the construction of confidence bands for an underlying cdf, which was developed by \cite{KanSri1972} and makes use of so called Kolmogorov-Smirnov type statistics. These statistics are of the form
\begin{equation*}
	K_{\check{\bvt}}\,=\,\sup_{x\in\R}|F_{\check{\bvt}}(x)-F_{\bvt}(x)|\,,
\end{equation*}
where $\check{\bvt}=\check{\bvt}(\bX)$ denotes an estimator of $\bvt$ based on $\bX$. If $\check{\bvt}$ is equivariant, i.e., if for all $a\boldsymbol{1}=(a,\dots,a)\in\R^m$ and $b>0$
\[\check{\bvt}(a\boldsymbol{1}+b\bX)=(a,0)+b\check{\bvt}(\bX)\,,\]
the statistic can be rewritten as
\begin{eqnarray}
	K_{\check{\bvt}}\,&=&\,
	\sup_{z\in\R} |F_{\check{\bvt}(\bX)}(\sigma z+\mu)-F_{(0,1)}(z)|\notag\\[1ex]
	\,&=&\,
	\sup_{z\in\R}|F_{\check{\bvt}((\bX-\mu\boldsymbol{1})/\sigma)}(z)-F_{(0,1)}(z)|\,.\label{eq:KSTpivot}
\end{eqnarray}
Here, $(\bX-\mu\boldsymbol{1})/\sigma$ is distributed as a progressively type-II censored sample with underlying cdf $F_{(0,1)}\in\mathcal{F}$, such that the distribution of $K_{\check{\bvt}}$ is found to be free of $\bvt$, i.e., $K_{\check{\bvt}}$ is a pivotal quantity.

The MLE $\hat{\bvt}$ of $\bvt$ is known to be equivariant; see, e.g., \cite{BalCra2014}, Section 12.1.1. This yields the following theorem, the proof of which is provided in the appendix; see also \cite{SriKanWha1975}, Section 2.1, for the complete sample case.

\begin{theorem}\label{thm:BandKST}
	Let $p\in(0,1)$. Then,
	\begin{equation}\label{eq:B_K}
		B_4\,=\,\left\{(x,y)\in\R^2:\,|F_{\hat{\bvt}}(x)-y|\leq d_{p}\right\}
	\end{equation}
	forms a confidence band for $F_{\bvt}$ with exact confidence level $1-p$, where $d_p\equiv d_p(m,n)$ denotes the $(1-p)$-quantile of
	\begin{equation}\label{eq:KSTrep}
		K_{\hat{\bvt}}\,\stackrel{\text{d}}{=}\,\max\{U,V\}
	\end{equation}
	with random variables
	\begin{eqnarray*}
		U\,&=&\,1-\exp\{-S\}\,,\\[1ex]
		V\,&=&\,|1-T|\,\exp\left\{\frac{S-T\ln(T)}{T-1}\right\}\,\left(\mathbbm{1}_{\{T<1\}}+\mathbbm{1}_{\{S<\ln(T)\}}\right)\,,
	\end{eqnarray*}
	where $S\thicksim F_{(0,1/n)}\in\mathcal{F}$ and $T\thicksim\Gamma(m-1,1/m)$ are independent. Here, $\stackrel{\text{d}}{=}$ means equality in distribution, and $\mathbbm{1}_A$ denotes the indicator function of the set $A$.
\end{theorem}

Note that Theorem \ref{thm:BandKST} allows for a simple numerical computation of the quantiles of $K_{\hat{\bvt}}$ by using Monte Carlo simulation. For some configurations, respective numerical values are presented in Table \ref{table:dp}.

\begin{table}[h!]
	\centering
	\begin{tabular}{cc|rrrrrrr}
		\toprule
		&& \multicolumn{7}{c}{$n$}\\[1ex] 
		&& 3 & 4 & 5 & 10 & 15 & 20 & 50 \\
		\midrule
		\multirow{7}{*}{$m$} & 3 & .123 & .109 & .099 & .075 & .064 & .058 & .045 \\ 
		& 4  &  & .095 & .086 & .064 & .055 & .049 & .037 \\ 
		& 5  &  &  & .078 & .058 & .049 & .044 & .033 \\ 
		& 10 &  &  &  & .045 & .038 & .034 & .024 \\ 
		& 15 &  &  &  &  & .033 & .029 & .020 \\ 
		& 20 &  &  &  &  &  & .027 & .018 \\ 
		& 50 &  &  &  &  &  &  & .014 \\ 
		\midrule
	\end{tabular}
	\caption{Value of $d_p\equiv d_p(m,n)$ for $B_4$ to have exact confidence level $1-p=90\%$ for different configurations $m,n$ (obtained by $10^7$ Monte Carlo simulations, each).}
	\label{table:dp}
\end{table}

$B_4$ in formula (\ref{eq:B_K}) has the same non-random vertical width $2d_{p}$ at every $x\in\R$, which implies that there will be points $(x,y)\in B_4$ that cannot be part of any graph of a cdf lying wholly inside the band. These points, however, may be removed from the band in a second step without affecting its confidence level. That is, if $B_4$ has exact confidence level $1-p$, the trimmed confidence band
\[B_4'\,=\,\bigcup_{\tilde{\bvt}\in\Theta:\,\text{graph}\,F_{\tilde{\bvt}}\subseteq B_4} \text{graph}\,F_{\tilde{\bvt}}\]
has also exact confidence level $1-p$. 

Note that $B_4'$ can be considered to be constructed from the confidence region 
\begin{align}
	C_4'
	&=\{\tilde{\bvt}\in\Theta:\,\text{graph}\,F_{\tilde{\bvt}}\subseteq B_4 \} \notag\\[1ex]
	&= \{ \tilde{\bvt}\in\Theta: \, \sup_{x\in\R}|F_{\tilde{\bvt}}(x)-F_{\hat{\bvt}}(x)| \leq d_p \}\label{eq:C4'}
\end{align}
in the sense of the first approach, i.e., $B_4'=B_{C_4'}$. Since $C_4'$ is exhaustive by definition, it follows that 
\begin{equation*}
	P_{\bvt}(\bvt \in C_4')\,=\,P_{\bvt}(\text{graph}\, F_{\bvt}\subseteq B_4') =1-p\,,\qquad\bvt\in\Theta\,,
\end{equation*}
such that $C_4'$ defines a confidence region for $\bvt$ with exact confidence level $1-p$. To compute $C_4'$, we may use that 
\begin{eqnarray*}
	\sup_{x\in\R}|F_{(\tilde{\mu},\tilde{\sigma})}(x) - F_{(\hat{\mu}, \hat{\sigma})}(x)|
	\,=\, \sup_{x\in\R}|F_{(\frac{\tilde{\mu}-\hat{\mu}}{\hat{\sigma}},\frac{\tilde{\sigma}}{\hat{\sigma}})}(x) - F_{(0,1)}(x)|\,, 
\end{eqnarray*}
and the latter quantity is calculated explicitly in the proof of Theorem \ref{thm:BandKST}.
The shape of $C_4'$ is depicted in Figure \ref{fig:C4}.

$C_4'$ can be trimmed even further without changing its confidence level by using that $\mu\leq\hat{\mu}$ $P_{\bvt}$-almost-surely. By defining
\begin{equation}\label{eq:C4''}
	C_4''\,=\,\{(\mu,\sigma)\in C_4':\mu\leq\hat{\mu}\}\quad\text{and}\quad B_4''\,=\,B_{C_4''}\,,
\end{equation}
we have for every $\bvt\in\Theta$ that
\begin{eqnarray*}
	1-p\,&=&\,P_{\bvt}(\bvt\in C_4')\,=\,P_{\bvt}(\bvt\in C_4'')\\[1ex]
	\,&\leq&\,P_{\bvt}(\text{graph}\,F_{\bvt}\subseteq B_4'')
	\,\leq\,P_{\bvt}(\text{graph}\,F_{\bvt}\subseteq B_4')\\[1ex]
	\,&=&\,1-p\,.
\end{eqnarray*}
Hence, the confidence region $C_4''$ and the corresponding confidence band $B_4''$ have both exact confidence level $1-p$.

$C_4'$ and $C_4''$ can be stated explicitly.

\begin{theorem}\label{thm:C4}
	The confidence regions $C'_4$ and $C''_4$ in formulas (\ref{eq:C4'}) and (\ref{eq:C4''}) admit the representations
	\begin{eqnarray*}
		C'_4 \,&=&\, \Big\{ (\mu,\sigma)\in \Theta:\, u(\tfrac{\hat{\sigma}}{\sigma}) \leq \frac{\hat{\mu}-\mu}{\sigma} \leq o(\tfrac{\hat{\sigma}}{\sigma})  \Big\}\,, \\[1ex]
		C''_4 \,&=&\, \Big\{ (\mu,\sigma)\in \Theta:\, \max\{u(\tfrac{\hat{\sigma}}{\sigma}),0\} \leq \frac{\hat{\mu}-\mu}{\sigma} \leq o(\tfrac{\hat{\sigma}}{\sigma})  \Big\}\,,
	\end{eqnarray*}
	where
	\begin{eqnarray*}
		u(x) \,&=&\, \begin{cases}
			\qquad h(x)\,, & x< 1-d_p\\
			x\ln(1-d_p),& x\geq 1-d_p
		\end{cases}\,, \\[1ex]
		o(x) \,&= &\,\begin{cases}
			-\ln(1-d_p),& x\leq\frac{1}{1-d_p} \\
			\qquad h(x)\,,& x>\frac{1}{1-d_p}
		\end{cases}\,,
	\end{eqnarray*}
	for $x>0$ and the function $h:(0,\infty)\rightarrow\R$ is defined by \begin{align*}
		h(x) = \ln\left( \frac{d_p}{|1-x|} \right) (x-1) + x\ln (x)\,,\quad x>0\,,x\neq1\,,
	\end{align*}
	and $h(1)=0$.
\end{theorem}

The confidence regions $C'_4$ and $C''_4$ are depicted in Figure \ref{fig:C4} along with the corresponding confidence bands $B'_4$, $B''_4$ and the untrimmed confidence band $B_4$. Here, the form of $B'_4$ and $B''_4$ is determined by maximizing/minimizing $F_{\bvt}(x)$ numerically on every confidence region for each $x\in\R$.

\section{Data Example}\label{s:comp}

For illustration and to compare the bands in terms of their width and area, we apply the confidence bands proposed in Section \ref{s:bands} to a progressively type-II censored data set presented by \cite{VivBal1994} and generated from a real data set in \cite{Nel1982}, pp. 105, 228. The observations are shown in Table  \ref{table:realdata} and consist of the times (in minutes) to breakdown of $m=8$ out of $n=19$ insulating fluids between electrodes at voltage 34 Kilovolts. The corresponding censoring scheme is given by $\boldsymbol{R}=(0,0,3,0,3,0,0,5)$. 

\begin{figure}[h!]
	\centering
	\includegraphics[width=0.49\textwidth]{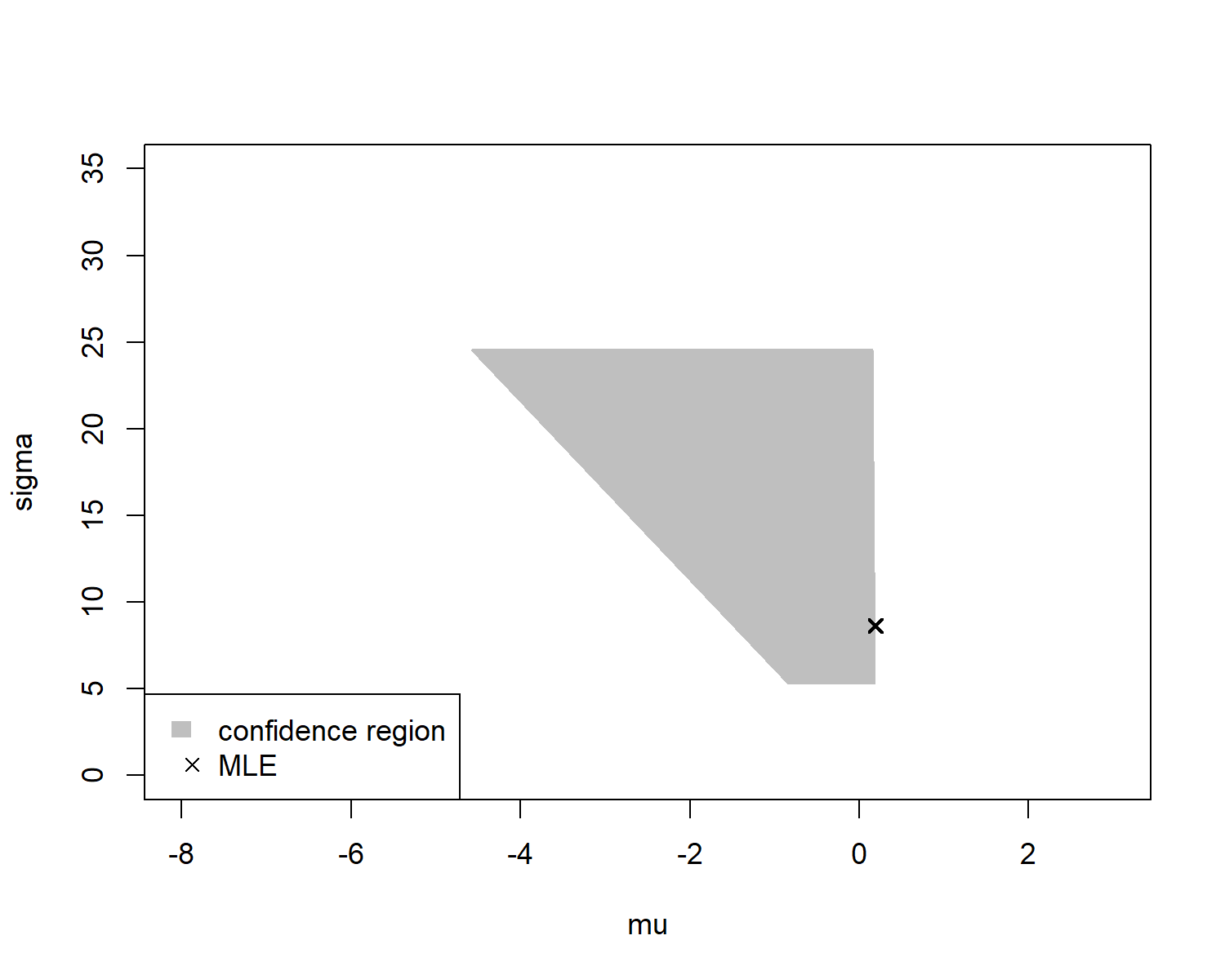}
	\includegraphics[width=0.49\textwidth]{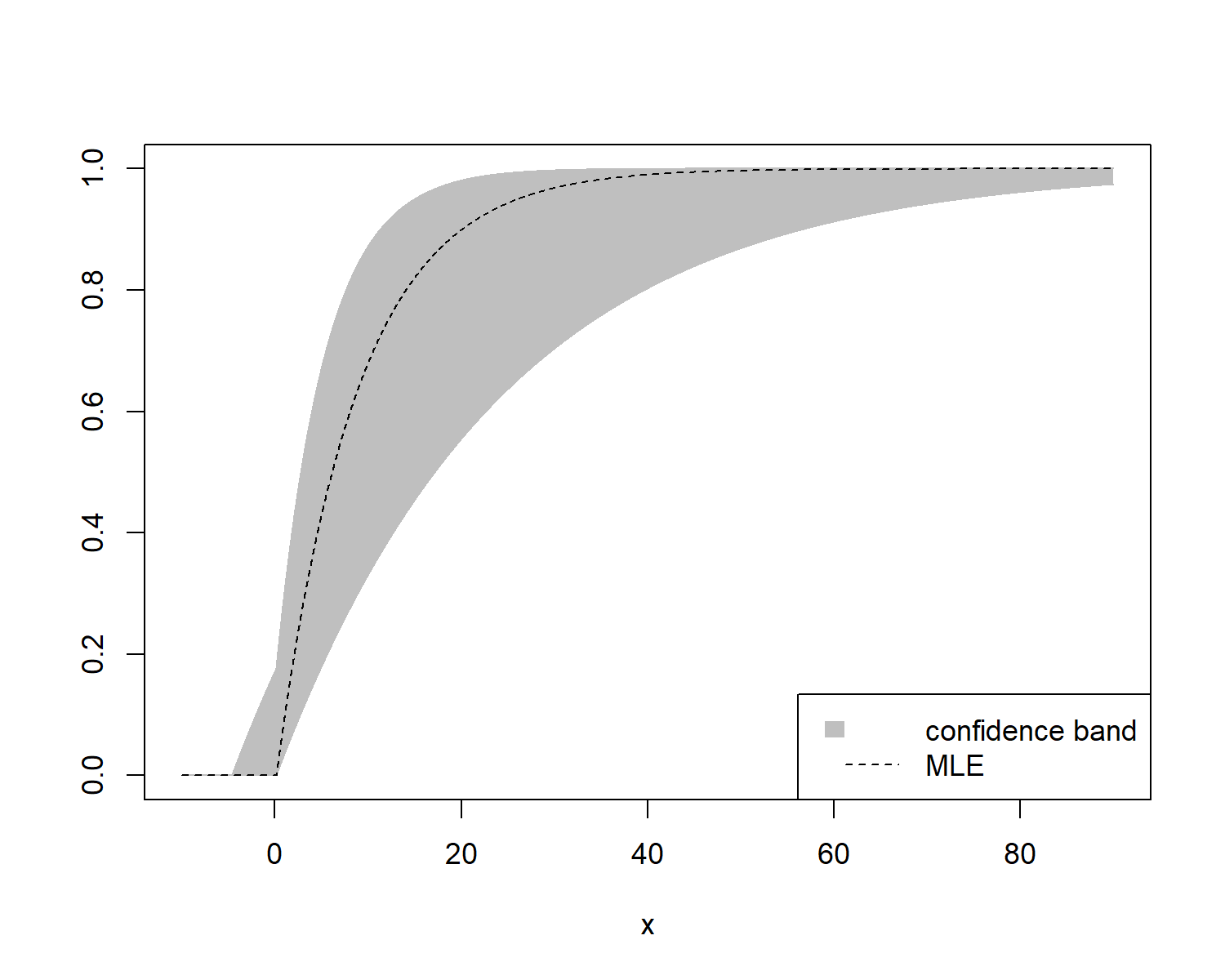}
	\caption{Confidence region $C_1$ and the resulting confidence band $B_1$ with exact confidence level 90.25\%, each, based on the insulating fluid data in Table \ref{table:realdata}.}
	\label{fig:C1}
\end{figure}

\begin{figure}[h!]
	\includegraphics[width=0.49\textwidth]{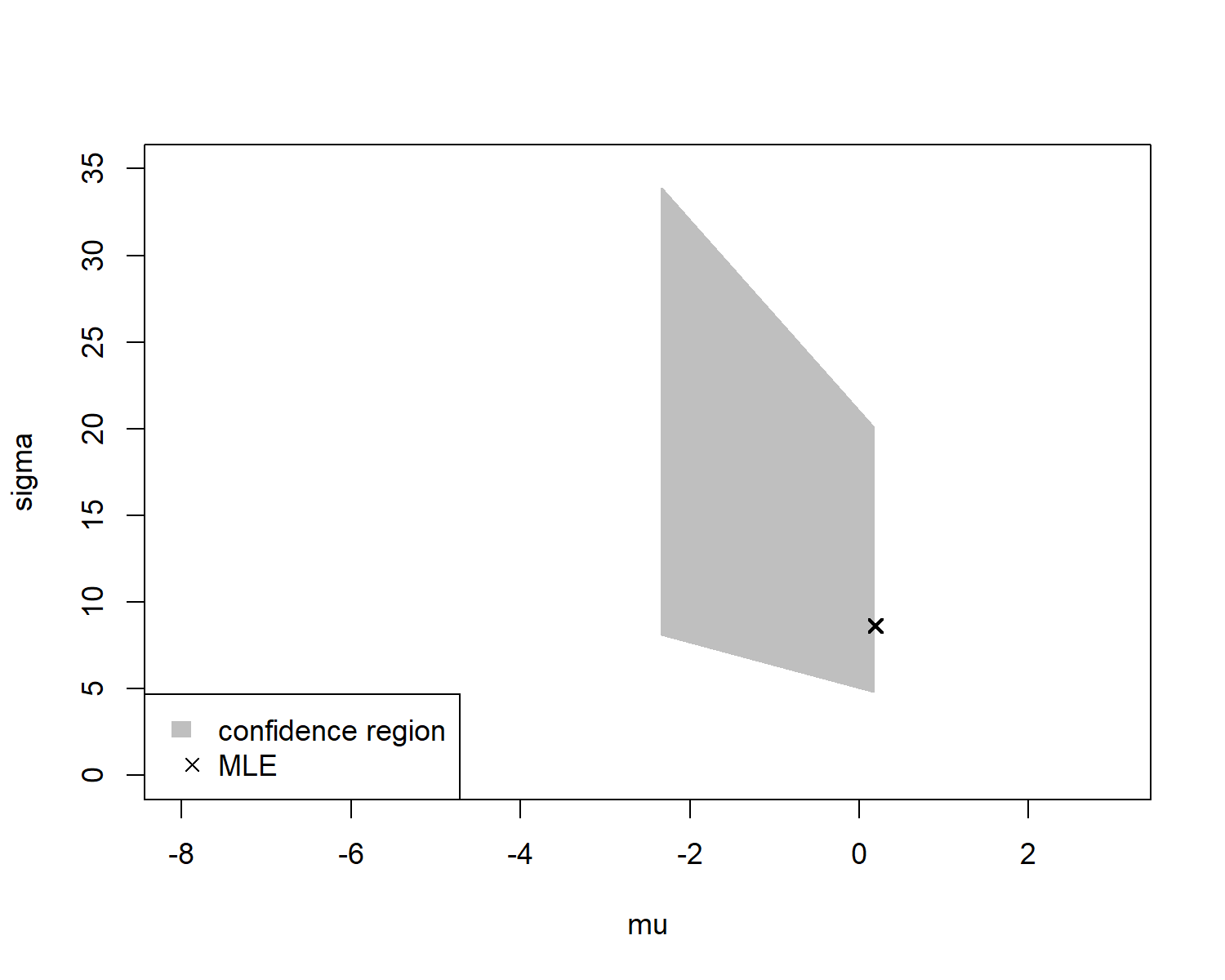}
	\includegraphics[width=0.49\textwidth]{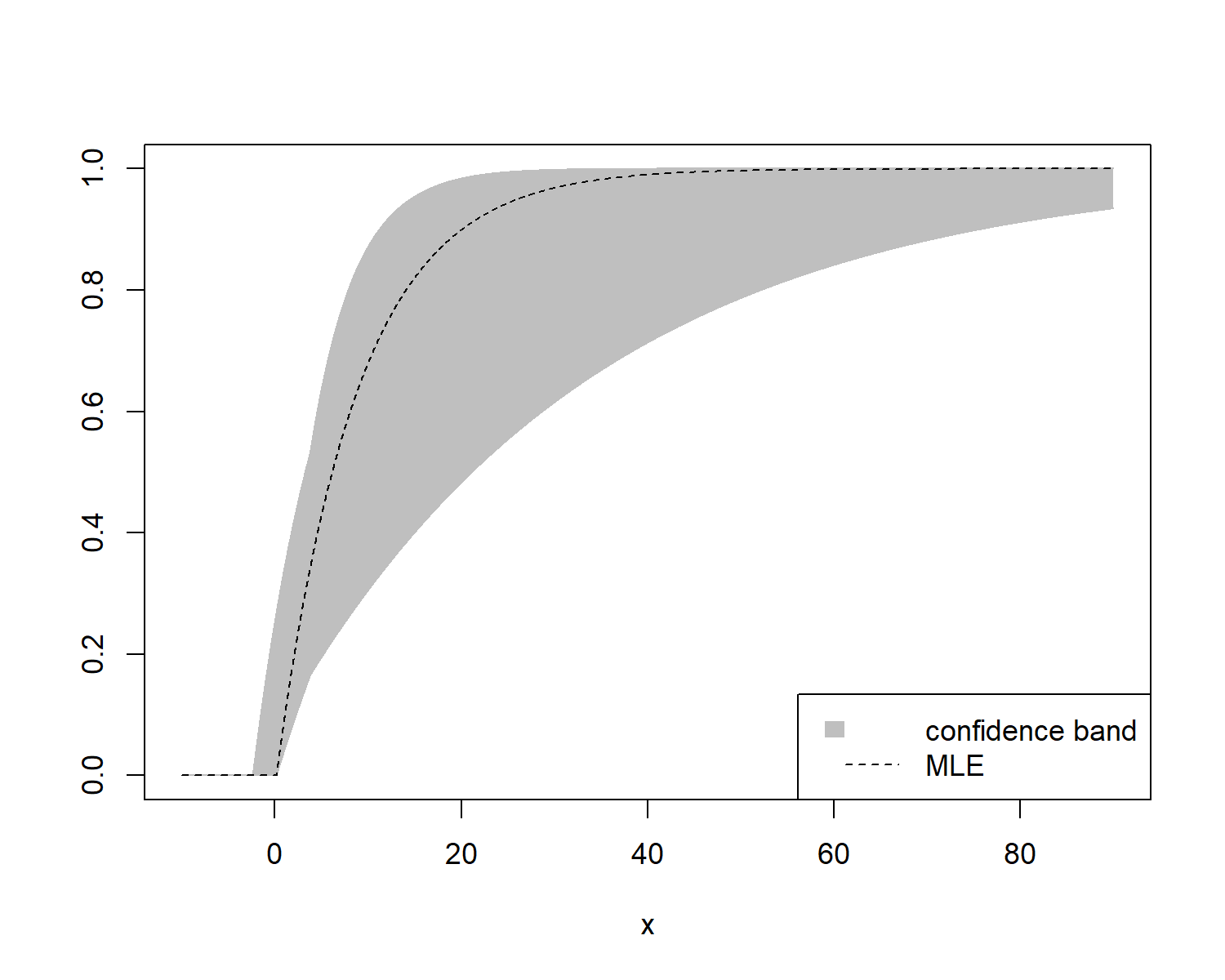}
	\caption{Confidence region $C_2$ and the resulting confidence band $B_2$ with exact confidence level 90.25\%, each, based on the insulating fluid data in Table \ref{table:realdata}.}
	\label{fig:C2}
\end{figure}

\begin{figure}[h!]
	\includegraphics[width=0.49\textwidth]{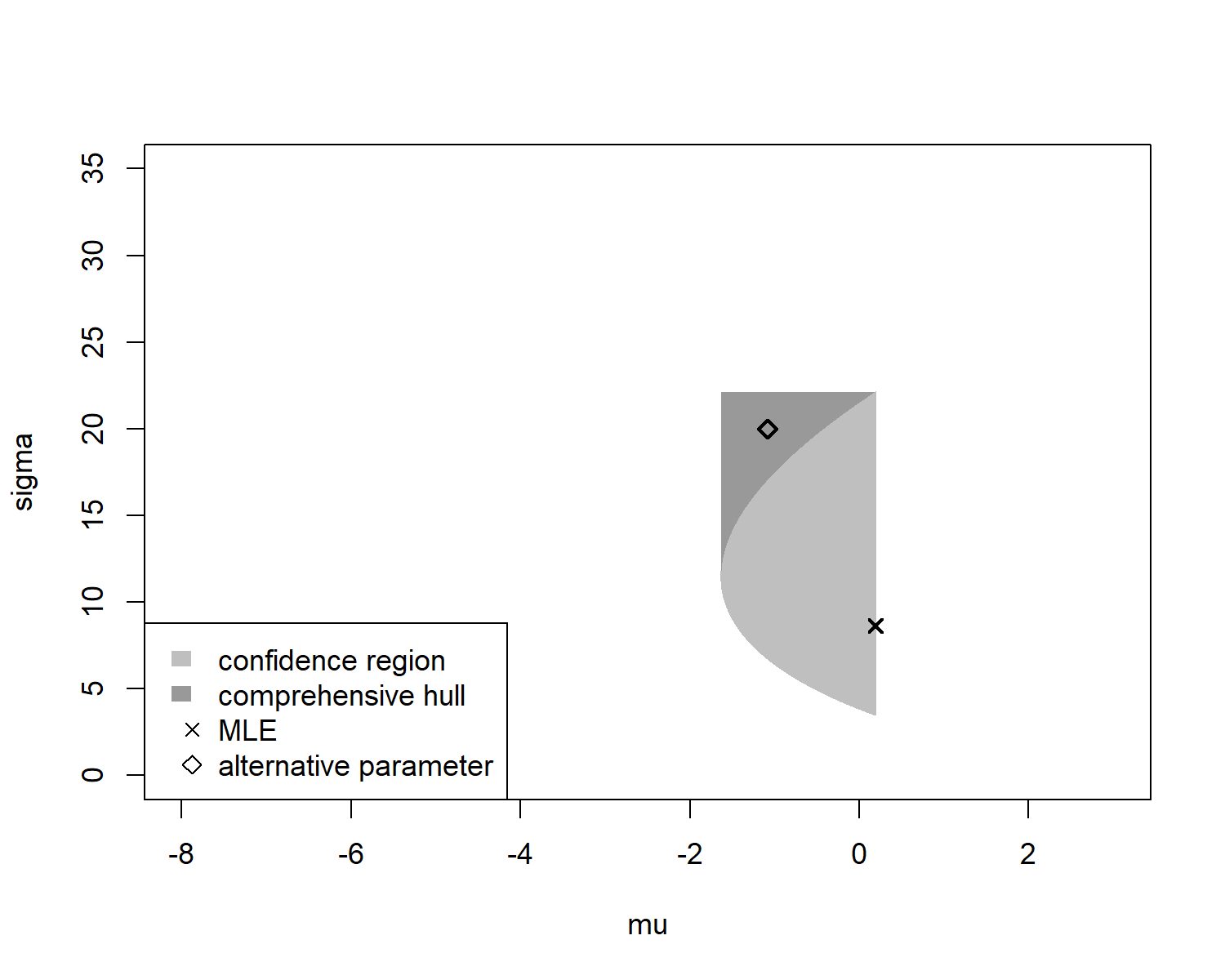}
	\includegraphics[width=0.49\textwidth]{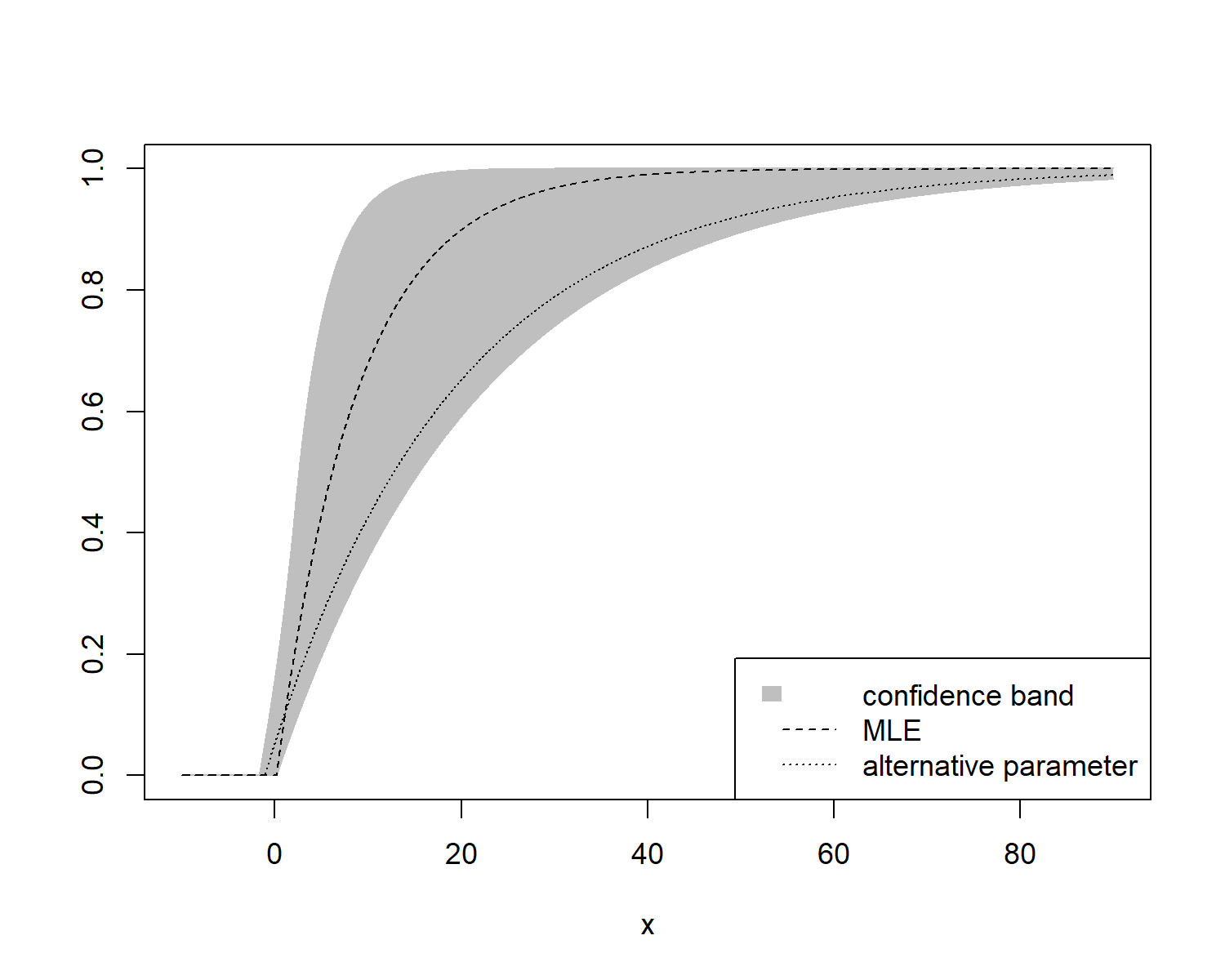}
	\caption{Confidence region $C_3$ with exact confidence level $1-p(\tau)=87.3\%$ and the resulting confidence band $B_3$ with exact confidence level $\tau=90.25\%$, based on the insulating fluid data in Table \ref{table:realdata}. $B_3$ covers the cdf associated with the marked alternative parameter, which is not contained in $C_3$ but part of its comprehensive convex hull.} 
	\label{fig:C3}
\end{figure}

\begin{figure}[h!]
	\includegraphics[width=0.49\textwidth]{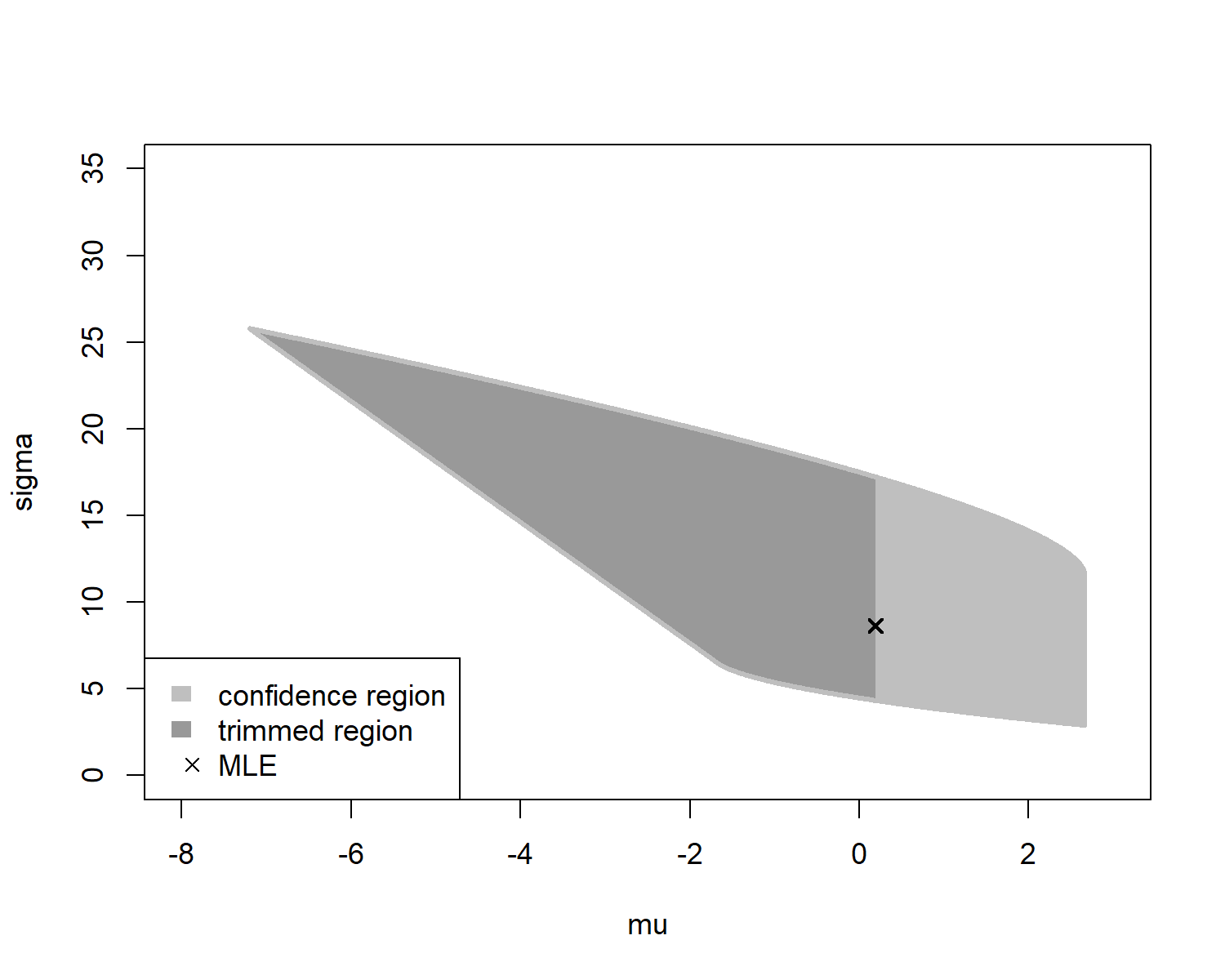}
	\includegraphics[width=0.49\textwidth]{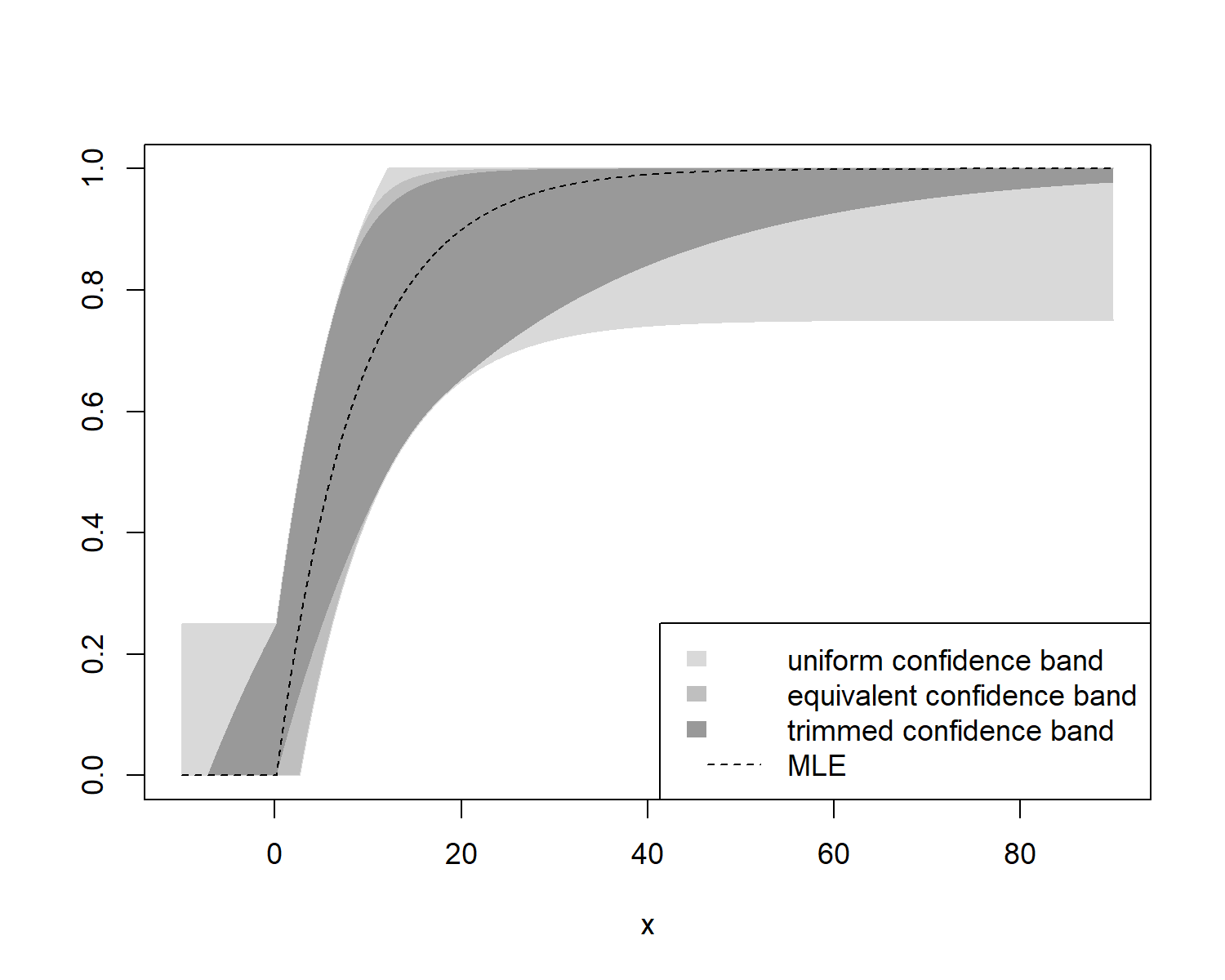}
	\caption{Uniform Kolmogorov-Smirnov type confidence band $B_4$ (right, light grey) and the resulting confidence regions $C_4'$ and $C_4''$ (left, grey and dark grey) with exact confidence level 90.25\% based on the insulating fluid data in Table \ref{table:realdata}. The corresponding confidence bands $B'_4$ and $B''_4$ are also depicted (right, grey and dark grey).} 
	\label{fig:C4}
\end{figure}

We assume that $x_{1:8:19},\dots,x_{8:8:19}$ in Table \ref{table:realdata} are realizations of progressively type-II censored order statistics from an exponential cdf $F_{\bvt}$ as in formula (\ref{eq:locscale}); see also \cite{VivBal1994} and  \cite{BalCra2014}, p. 255. According to formulas (\ref{eq:mudach}) and (\ref{eq:sigmadach}), the MLEs of $\mu$ and $\sigma$ based on the data are given by $\hat{\mu}=.19$ and $\hat{\sigma}=8.635$. For an exact confidence level of 90.25\% $(=95\%\times95\%)$, the boundaries of the confidence bands $B_1,B_2,B_3$, $B_4$, $B_4'$, and $B_4''$ are depicted in Figures \ref{fig:C1}-\ref{fig:C4}. Here, for the underlying confidence regions $C_1$ and $C_2$,  the exact confidence level is uniformly allocated to the intervals and tails, i.e., we set $q_1=2.5\%$ and $q_2=97.5\%$ in formulas (\ref{eq:Wu1}) and (\ref{eq:Wu2}), respectively. Moreover, to ensure that $B_3$ meets the exact confidence level of $\tau=90.25\%$, we choose the constant $c_p$  in formula (\ref{eq:Cmin}) as $c_{p(\tau)}$, which is determined by using Lemma \ref{La:CPMin} yielding $1-p(\tau)=87.3\%$ and $c_{p(\tau)}=-11.587$. The constant $d_{0.9025}=0.249$ in formula (\ref{eq:B_K}), in turn, is numerically obtained by sampling according to Theorem \ref{thm:BandKST} with $10^7$ simulations.


\begin{table*}[h!]
	\centering		
	\begin{tabular}{cc|rrrrrrrr}
		\toprule
		&$i$ & 1 & 2 & 3 & 4 & 5 & 6 & 7 & 8\\ \midrule
		&$x_{i:8:19}$ & .19 & .78 & .96 & 1.31  & 2.78 & 4.85 & 6.50 & 7.35\\[1ex]
		&$R_i$ & 0 & 0 & 3 & 0 & 3 & 0 & 0 & 5\\
		\midrule
	\end{tabular}		
	\caption{Progressively type-II censored data set in \cite{VivBal1994}: times (in minutes) to breakdown of $m=8$ out of $n=19$ insulating fluids at 34 Kilovolts with censoring scheme $\boldsymbol{R}=(R_1,\dots,R_8)$.}	
	\label{table:realdata}
\end{table*}

Numerical values for the maximum band width $W(B)$ and the area $A(B)$ of $B_1,\dots,B_4$, $B_4'$, and $B_4''$ are shown in Table \ref{table:realdataprops}. It gives the rankings
\begin{eqnarray*}
	W(B_4'')\leq W(B_4')\leq W(B_4) \leq W(B_1) \leq W(B_2) \leq W(B_3)
\end{eqnarray*}
and
\begin{eqnarray*}
	A(B_4'')\leq A(B_4')\leq A(B_3)\leq A(B_1) \leq A(B_2) \leq A(B_4)\,.
\end{eqnarray*}
In particular, the trimmed Kolmogorov-Smirnov type bands $B_4''$ and $B_4'$ are first and second best for both criteria; note that the untrimmed band $B_4$ has infinite area by construction. While the area of $B_3$, which is obtained from a minimum area confidence region for the underlying parameter, is also comparatively small, its maximum band width is worst among the confidence bands considered. Regarding the confidence bands based on the trapezoidal confidence regions, $B_1$ performs better than $B_2$.

\begin{table}[h!]
	\centering		
	\begin{tabular}{cc|rrrrrr}
		\toprule
		&$B$ & $B_1$ & $B_2$ & $B_3$ & $B_4$ & $B_4'$ & $B_4''$ \\ \midrule
		&$W(B)$ & .54 & .57 & .59 & .50 & .50 & .47\\
		&$A(B) $ & 20.59& 27.53 & 18.87 & $\infty$ & 18.70 & 17.90\\ 			\midrule
	\end{tabular}		
	\caption{Maximum band width $W(B)$ and area $A(B)$ of the confidence bands $B_1,\dots,B_4$, $B_4'$, and $B_4''$ with exact confidence level 90.25\% based on the insulating fluid data in Table \ref{table:realdata}.}	
	\label{table:realdataprops}
\end{table}


\section{Generalizations and Extensions}\label{s:general}

Generalizations and extensions of the preceding results are possible in different directions, which are pointed out in the following. First, we focus on related models for ordered data, in which the findings may be applied as well with minor changes. Then, we suppose to have an underlying location-scale family of distributions containing the exponential one as particular case.

\subsection{Sequential Order Statistics}\label{ss:sos}
In the above derivations, the explicit structure of $\gamma_1,\dots,\gamma_m$ defined by formula (\ref{eq:gammas}) via the censoring scheme has not been used. Indeed, all findings remain true as long as the $\gamma$'s are given positive parameters. Here, the particular choice
\begin{equation}\label{eq:gammas2}
	\gamma_j\,=\,(n-j+1)\,\alpha_j\,,\qquad 1\leq j\leq n\,,
\end{equation}
with known positive parameters $\alpha_1,\dots,\alpha_n$ is of some interest. In that case, $\bX$ is distributed as the first $m$ (of $n$) sequential order statistics (SOSs) based on the cdfs
\begin{eqnarray*}
	F_j(x)\,&=&\,1-(1-F_{\bvt}(x))^{\alpha_j}\\[1ex]
	\,&=&\,1-\exp\left\{-\frac{\alpha_j(x-\mu)}{\sigma}\right\}\,,\qquad x>\mu\,,
\end{eqnarray*}
with corresponding hazard rates $\lambda_{F_j}=\alpha_j\lambda_{F_{\bvt}}=\alpha_j/\sigma$ on $(\mu,\infty)$, $1\leq j\leq n$, being proportional to that of $F_{\bvt}$; see \cite{Kam1995a,Kam1995b}. The model of SOSs allows for describing component lifetimes of sequential $k$-out-of-$n$ systems, which are operating as long as $k$ out of $n$ components are operating, where the failure of some component may have an impact on the residual lifetimes of the remaining components ($k=n-m+1$, here). In particular, they may be used to model load-sharing effects arising in systems, in which the components share some total load. All $n$ components then start working at hazard rate $\alpha_1\lambda_{F_{\bvt}}$ and, upon the $j$-th failure, $1\leq j\leq m-1$, the hazard rate of the surviving components changes to $\alpha_{j+1}\lambda_{F_{\bvt}}$. The $m$-th SOS $X_m$ coincides with the system failure time, after which no further observations are recorded such that the data is type-II right censored. Note that order statistics (based on $F_{\bvt}$) modelling the component lifetimes of the common $k$-out-of-$n$ system are contained in the model by setting $\alpha_1=\dots=\alpha_n\,(=1)$. In this context, confidence bands for the baseline cdf $F_{\bvt}$ may be useful to assess the quality of the individual components.  

Having said that the $m$-th SOS $X_m$ also represents the lifetime of the system, a confidence band for its cdf is naturally of interest as well. For arbitrary choices of $\alpha_1,\dots,\alpha_n$, closed formula expressions of the marginal cdf $F_{\bvt}^{(m)}$ of $X_m$
are derived in \cite{CraKam2003} by using Meijer's $G$-functions. If $\alpha_1,\dots,\alpha_n$ are such that the $\gamma$'s in formula (\ref{eq:gammas2}) are pairwise distinct,
this formula simplifies to $F_{\bvt}^{(m)}=H(F_{\bvt})$ with function
\begin{equation*}
	H(y)\,=\, 1-\left( \prod_{j=1}^m \gamma_j \right) \sum_{i=1}^m \frac{a_i}{\gamma_i} \left( 1-y \right)^{\gamma_i}\,,\qquad y\in[0,1]\,,
\end{equation*}
where
\[a_i \,=\, \prod_{j=1, j\neq i}^m \frac{1}{\gamma_j-\gamma_i}\,,\qquad 1\leq i\leq m\,;
\]
see, e.g., \cite{CraKam2001b}, Theorem 2.5. 
In fact, $H$ coincides with the cdf on $[0,1]$ of the $m$-th SOS based on a standard uniform distribution and model parameters $\alpha_1,\dots,\alpha_n$ (i.e., based on the cdfs $F_j(x)=1-(1-u)^{\alpha_j}$, $u\in[0,1]$, $1\leq j\leq n$); $H:[0,1]\rightarrow[0,1]$ is therefore strictly increasing and, in particular, bijective. This transformation can be used for the construction of a confidence band for $F_{\bvt}^{(m)}$ as follows.

Suppose that we have already constructed a confidence band $B$, say, for $F_{\bvt}$ with (exact) confidence level $1-p\in(0,1)$. Let
\[\breve{B}\,=\,\{(x,y)\in\R\times[0,1]:\,(x,H^{-1}(y))\in B\}\,.\]
Then,
\begin{eqnarray*}
	\text{graph}\, F_{\bvt}^{(m)} \subseteq \breve{B}\quad&\Leftrightarrow&\quad
	(x,H(F_{\bvt}(x)))\in\breve{B}\;\;\forall x\in\R\\[1ex]
	\quad&\Leftrightarrow&\quad
	(x,F_{\bvt}(x))\in B\;\;\forall x\in\R\\[1ex]
	\quad&\Leftrightarrow&\quad
	\text{graph}\, F_{\bvt}\subseteq  B\,.
\end{eqnarray*}
Hence, $\breve{B}$ forms a confidence band for $F_{\bvt}^{(m)}$ with (exact) confidence level $1-p$.

More results on inference with SOSs can be found, for instance, in \cite{CraKam2001b} and \cite{Cra2016}. For recent works, see, e.g., \cite{BurCraGor2016} and \cite{MieBed2019}. Finally, notice that SOSs with proportional hazard rates are closely related to generalized order statistics; see \cite{Kam1995a,Kam1995b,Kam2016}. Some of the results presented here may therefore also be interpreted in light of other models of ordered random variables and may be applied, for example, to construct a confidence band for the exponential baseline cdf of Pfeifer record values.

\subsection{A Location-Scale Familiy of Distributions}\label{ss:otherlsf}

In the context of SOSs, various inferential results are shown in \cite{CraKam2001b} for the underlying cdf belonging to the location-scale family
\begin{equation}\label{eq:locscale2}
	F_{\bvt}(x)\,=\,1-\exp\left\{-\,\frac{g(x)-\mu}{\sigma}\right\}\,,\quad x\geq g^{-1}(\mu)\,,
\end{equation}
for $\bvt=(\mu,\sigma)\in (\lim_{u\downarrow0}g(u),\infty)\times(0,\infty)$, where $g$ is a differentiable and strictly increasing function on $(0,\infty)$ satisfying $\lim_{x\rightarrow\infty}g(x)=\infty$. Note that the supports of the distributions in that family are bounded from the left. Exponential distributions and Pareto distributons are included in this setting and are obtained by choosing $g(x)=x$ and $g(x)=\ln(x)$, respectively; more examples can be found in \cite{CraKam2001b}, Section 6. When replacing assumption (\ref{eq:locscale}) by (\ref{eq:locscale2}) with a known function $g$, the MLEs of $\mu$ and $\sigma$ in formulas (\ref{eq:mudach}) and (\ref{eq:sigmadach}) change to
\begin{eqnarray*}
	\hat{\mu}\,&=&\ g(X_{1:m:n})\\[1ex]
	\text{and}\quad\hat{\sigma}\,&=&\,\frac{1}{m}
	\sum_{j=2}^m \gamma_j(g(X_{j:m:n})-g(X_{j-1:m:n}))\,,
\end{eqnarray*}
where all distributional properties of $\hat{\mu}$ and $\hat{\sigma}$ are preserved; see, e.g., \cite{CraKam2001b}, Section 9. Since the construction principles and arguments in Section \ref{s:bands} do not use the explicit form of the MLEs, all findings remain true under this setup.


\section{Conclusion}\label{s:conclusion}

Based on a progressively type-II censored sample from the two-parameter exponential distribution, various confidence bands are derived containing the entire graph of the underlying cumulative distribution function with a desired exact probability. 
The bands are constructed via confidence regions for the location-scale parameter and Kolmogorov-Smirnov type statistics, and the relation between both approaches is highlighted. 
Explicit formulas for the boundaries of the confidence bands as well as instructions how to obtain the quantiles of the relevant statistics by simulation are provided, such that the bands may easily be computed in applications. The Kolmogorov-Smirnov type statistic can also be used to construct a novel confidence region for the location-scale parameter, and an explicit representation of the latter is found. By means of a data example, the bands are illustrated and compared in terms of maximum band width and area, where the trimmed Kolmogorov-Smirnov type band is found to perform well for both criteria. 
The presented results also yield useful findings for related models of ordered random variables. 
This comprises, for instance, confidence bands for the baseline and marginal cumulative distribution functions of sequential order statistics, which serve as a model for the component lifetimes of sequential $k$-out-of-$n$ systems.
Finally, the results are extended to some other location-scale families of distributions.


\section*{Appendix}
\subsection{Proof of Theorem \ref{thm:Wu}}
To find the boundaries of $B_2$, we proceed similarly as in \cite{SriKanWha1975}, Section 2.2. Let $x\in\R$ be fixed. The aim is to determine the parameters $\bvt\in C_2$, which maximize and minimize $F_{\bvt}(x)$. This question may be simplified to finding $(\mu,\sigma)\in C_2$ corresponding to extremal values of $(x-\mu)/\sigma$ . Since the latter is a monotone function of $\sigma$ for fixed $\mu$, these points lie on the graphs of the functions $\sigma=b_{q_i}(\mu)$, $\mu\in[\mu_{q_2},\mu_{q_1}]$, $i=1,2$. There we have that
\begin{equation*}
	\frac{x-\mu}{\sigma}\,=\,\frac{(x-\mu)\chi_{q_i}^2(2m)}{2(n(\hat{\mu}-\mu)+m\hat{\sigma})}
\end{equation*}
is strictly decreasing in $\mu$ if $x\leq\hat{\mu}+m\hat{\sigma}/n$, and strictly increasing in $\mu$ if $x>\hat{\mu}+m\hat{\sigma}/n$. Hence,
\begin{eqnarray*}
	\operatorname*{arg\,min}_{\bvt\in C_2}\frac{x-\mu}{\sigma}
	\,=\,\begin{cases} (\mu_{q_1},b_{q_2}(\mu_{q_1}))\,, & x\leq \mu_{q_1}\\[1ex]
		(\mu_{q_1},b_{q_1}(\mu_{q_1}))\,, & \mu_{q_1}<x\leq \hat{\mu}+m\hat{\sigma}/n\\[1ex]
		(\mu_{q_2},b_{q_1}(\mu_{q_2}))\,, & x>\hat{\mu}+m\hat{\sigma}/n\end{cases}\,,
\end{eqnarray*}
and
\begin{eqnarray*}
	\operatorname*{arg\,max}_{\bvt \in C_2}\frac{x-\mu}{\sigma}\,=\,\begin{cases} (\mu_{q_2},b_{q_1}(\mu_{q_2}))\,, & x\leq \mu_{q_2}\\[1ex]
		(\mu_{q_2},b_{q_2}(\mu_{q_2}))\,, & \mu_{q_2}<x\leq \hat{\mu}+m\hat{\sigma}/n\\[1ex]
		(\mu_{q_1},b_{q_2}(\mu_{q_1}))\,, & x>\hat{\mu}+m\hat{\sigma}/n\end{cases}\,.
\end{eqnarray*} 
This yields the representations of $U_2$ and $O_2$, where each the first two cases can be summarized, since $F_{(\mu_{q_i},\sigma)}(x)=0$ for $x\leq\mu_{q_i}$, $\sigma>0$, and $i=1,2$. 

\subsection{Proof of Theorem \ref{thm:minvol}}

Let $x\in\R$ be fixed. As in the proof of Theorem \ref{thm:Wu}, we have to find the extremal values of $(x-\mu)/\sigma$ with respect to $(\mu,\sigma)$ lying on the boundary of $C_3$, which is given by the vertical line $\hat{\mu}\times[Z_{-1},Z_0]$ and the curve $\{(\hat{\mu}+g(\sigma),\sigma):\sigma\in[Z_{-1},Z_0]\}$. Clearly, on the line, $(x-\mu)/\sigma=(x-\hat{\mu})/\sigma$ is extremal only at $\sigma=Z_{-1}$ and $\sigma=Z_0$. On the curve, we have
\begin{eqnarray*}
	\frac{x-\mu}{\sigma}\,=\,\frac{x-\hat{\mu}}{\sigma}-\frac{(m+1)\ln(\sigma)}{n}-\frac{m\hat{\sigma}}{n\sigma}+d
\end{eqnarray*}
for some constant $d$ free of $\mu$ and $\sigma$, and, considered a function of $\sigma>0$, the right-hand side is strictly increasing on $(0,\sigma_x^*]$ and strictly decreasing on $[\sigma_x^*,\infty)$. Together, these findings imply that, in any case, the minimum of $(x-\mu)/\sigma$ over the boundary of $C_3$ is attained on the line, and the respective maximum is attained on the curve. The case distinctions $x\leq\hat{\mu}$ and $x>\hat{\mu}$ for the minimum, and $\sigma_x^*<Z_{-1}$, $\sigma_x^*\in[Z_{-1},Z_0]$, and $\sigma_x^*>Z_0$ for the maximum then lead to the stated representation of $B_3$, where the first can eventually be dropped again, since $F_{(\hat{\mu},\sigma)}(x)=0$ for $x\leq\hat{\mu}$ and $\sigma>0$.

\subsection{Proof of Lemma \ref{La:CPMin}}
For $z\in[Z_{-1},Z_0]$, the function $g$ in formula (\ref{eq:g}) is strictly decreasing-increasing with minimum at $Z_{\text{min}}=\hat{\sigma}\exp\{-1-c_p/(m+1)\}$ given by
\begin{equation*}
	g(Z_{\text{min}})\,=\,\hat{\sigma}\left(\frac{m}{n}-\frac{m+1}{n}\exp\left\{-1-\frac{c_p}{m+1}\right\}\right)\,;
\end{equation*}
cf. \cite{LenBedKam2019}. Hence, the comprehensive convex hull of $C_3$ is given by the disjoint union $C_3\cup\Delta$, where
\begin{eqnarray*}
	\Delta\,&=&\,\left\{(\mu,\sigma)\in\Theta:\hat{\mu}+g(Z_{\text{min}})\leq\mu<\hat{\mu}+g(\sigma)\,,\, Z_{\text{min}}\leq\sigma\leq Z_0\right\}\\[1ex]
	\,&=&\,\left\{(\mu,\sigma)\in\Theta: (m+1)\ln\left(\frac{1}{m}\frac{m\hat{\sigma}}{\sigma}\right)-\frac{m\hat{\sigma}}{\sigma}-c_p\right.\\[1ex]
	&&\qquad \left.<\frac{n(\hat{\mu}-\mu)}{\sigma}
	\leq \left(\frac{m+1}{z}-1\right)\frac{m\hat{\sigma}}{\sigma}\;,\; y\leq \frac{m\hat{\sigma}}{\sigma}\leq z\right\}
\end{eqnarray*}
with $y$ and $z$ as stated (the set $\Delta$ is colored dark in Figure \ref{fig:C3}). Now, recall that the statistics $m\hat\sigma/\sigma\sim \Gamma(m-1,1)$ and $n(\hat\mu-\mu)/\sigma\sim \Gamma(1,1)$ are independent; see Section \ref{s:PCII}. Denoting by $g_k$ the density function of $\Gamma(k,1)$, we then obtain for $(\mu,\sigma)\in\Theta$
\begin{eqnarray*}
	P_{(\mu,\sigma)}((\mu,\sigma)\in\Delta)&&\,=\,\int \mathbbm{1}_\Delta\,dP_{(\mu,\sigma)}\\[1ex]
	\,=\,&&\int_y^z\int_{(m+1)\ln(v/m)-v-c_p}^{[(m+1)/z-1]v}
	g_1(u)g_{m-1}(v)\,du\,dv\\[1ex]
	\,=\,&&\frac{e^{c_p}m^{m+1}}{(m-2)!}\int_y^z v^{-3}\,dv\,-\,\frac{1}{(m-2)!}\int_y^z v^{m-2}e^{-(m+1)v/z}\,dv\,.
\end{eqnarray*}
A change of variables in the second integral then leads to the stated formula.

\subsection{Proof of Theorem \ref{thm:BandKST}}

By definition of $B_4$ and $d_p$, it is clear that for every $\bvt\in\Theta$
\begin{equation*}
	P_{\bvt}(\text{graph}\, F_{\bvt}\subseteq B_4)\,=\,P_{\bvt}(K_{\hat{\bvt}}\leq d_p)\,=\,1-p\,.
\end{equation*}
Hence, we only have to verify formula (\ref{eq:KSTrep}). To this end, let $\bvt=(\mu,\sigma)\in\Theta$ and
\begin{eqnarray*}
	K_{\bvt}\,&=&\,\sup_{x\in\R} |F_{\bvt}(x)-F_{(0,1)}(x)|\\[1ex]
	&=&\,\sup_{x>\min\{\mu,0\}} |F_{\bvt}(x)-F_{(0,1)}(x)|\,.
\end{eqnarray*}
A case distinction on the sign of $\mu$ and
resolving the indicator functions gives
\begin{eqnarray}
	K_{\bvt}\,&=&\,\max\left\{1-\exp\left\{-\max\left\{\mu,-\,\frac{\mu}{\sigma}\right\}\right\}\,,\,\sup_{x>\max\{\mu,0\}} |\kappa(x)|\right\}\label{eq:Kmax}
\end{eqnarray}
with function
\begin{equation*}
	\kappa(x)\,=\,\exp\{-x\}-\exp\left\{-\,\frac{x-\mu}{\sigma}\right\}\,,\qquad x\in\R.
\end{equation*}
Since 
\begin{equation*}
	\lim_{x\rightarrow\max\{\mu,0\}}|\kappa(x)|\,=\,1-\exp\left\{-\max\left\{\mu,-\,\frac{\mu}{\sigma}\right\}\right\}\
\end{equation*}
and $\lim_{x\rightarrow\infty} |\kappa(x)|=0$, the aim is now to find all local extrema of $\kappa$ in $(\max\{\mu,0\},\infty)$. For $\sigma=1$ and $\mu\neq0$, such extrema do not exist, since $\kappa$ is then strictly monotone. Hence, we have $K_{\bvt}=1-\exp\{-\max\{\mu,-\mu/\sigma\}\}$ for $\sigma=1$. For $\sigma\neq1$, simple analysis shows that $\kappa$ is either increasing/decreasing or decreasing/increasing with a local extrema at
\begin{equation*}
	x^*\,=\,\frac{\mu-\sigma\ln(\sigma)}{1-\sigma}\,.
\end{equation*}
The value 
\begin{equation*}
	|\kappa(x^*)|\,=\,|1-\sigma|\,\exp\left\{\frac{\mu-\sigma\ln(\sigma)}{\sigma-1}\right\} 
\end{equation*}
thus has to be taken into account in the calculation of the maximum in formula (\ref{eq:Kmax}) if
\begin{eqnarray*}
	x^*>\max\{\mu,0\}\qquad
	&\Leftrightarrow&\qquad
	\min\{x^*-\mu,x^*\}>0\\[1ex]
	\qquad&\Leftrightarrow&\qquad\min\left\{\,\frac{\mu-\ln(\sigma)}{1-\sigma}\,,\,\frac{\mu-\sigma\ln(\sigma)}{1-\sigma}\right\}>0\\[1ex]
	\qquad&\Leftrightarrow&\qquad\frac{\mu}{\ln(\sigma)}<\min\{\sigma,1\}\,.
\end{eqnarray*}
Hence, $K_{\bvt} = \max\{U(\bvt), V(\bvt)\}$ for $\sigma\neq1$, where
\begin{eqnarray*}
	U(\bvt)\,&=&\,1-\exp\left\{-\max\left\{\mu,-\,\frac{\mu}{\sigma}\right\}\right\}\,,\\[1ex]
	V(\bvt)\,&=&\,|1-\sigma|\,\exp\left\{\frac{\mu-\sigma\ln(\sigma)}{\sigma-1}\right\}\mathbbm{1}_{\left\{\frac{\mu}{\ln(\sigma)}<\min\{\sigma,1\}\right\}}\,.
\end{eqnarray*}
By setting $V(\bvt)=0$ for $\sigma=1$, the representation $K_{\bvt}=\max\{U(\bvt),V(\bvt)\}$ then  holds true for all $\bvt\in\Theta$.

Now, according to formula (\ref{eq:KSTpivot}), replacing $\bvt$ by 
\begin{equation*}
	\hat{\bvt}\left(\frac{\bX-\mu\boldsymbol{1}}{\sigma}\right)\,=\,\left(\frac{\hat{\mu}(\bX)-\mu}{\sigma},\frac{\hat{\sigma}(\bX)}{\sigma}\right)\,\stackrel{d}{=}\,(S,T)
\end{equation*}
completes the proof, where the distributional properties of $S$ and $T$ follow from formulas (\ref{eq:mudach}) and (\ref{eq:sigmadach}). Here, the final representations of $U$ and $V$ are obtained by using that $S$ and $T$ are positive almost-surely.

\subsection{Proof of Theorem \ref{thm:C4}}
For $\bvt=(\mu,\sigma)\in\Theta$, we have from the proof of Theorem \ref{thm:BandKST} that
\begin{equation}\label{eq:Kmax}
	K_{\bvt}\leq d_p\quad\Leftrightarrow\quad \max\{U(\bvt),V(\bvt)\}\leq d_p\,.
\end{equation}
First, note that
\begin{equation}\label{eq:U}
	U(\bvt)\leq d_p\quad\Leftrightarrow\quad 
	\sigma\ln(1-d_p)\leq\mu\leq-\ln(1-d_p)\,.
\end{equation}
Moreover, we introduce the function $h:(0,\infty)\rightarrow\R$ via \begin{align*}
	h(x) = \ln\left( \frac{d_p}{|1-x|} \right) (x-1) + x\ln (x)\,,\quad x>0\,,x\neq1\,,
\end{align*}
and $h(1)=0$, such that for $\sigma\leq1$
\begin{equation}\label{eq:V1}
	V(\bvt)\leq d_p\quad\Leftrightarrow\quad \mu\,\mathbbm{1}_{\{\mu>\sigma\ln(\sigma)\}}\geq h(\sigma)\,\mathbbm{1}_{\{\mu>\sigma\ln(\sigma)\}}\,,
\end{equation}
and for $\sigma>1$
\begin{equation}\label{eq:V2}
	V(\bvt)\leq d_p\quad\Leftrightarrow\quad \mu\,\mathbbm{1}_{\{\mu<\ln(\sigma)\}}\leq h(\sigma)\,\mathbbm{1}_{\{\mu<\ln(\sigma)\}}\,,
\end{equation}

Let us first find the lower bound $u(\sigma)$ for $\mu$ in dependence of $\sigma$, such that  the right-hand side of equivalence (\ref{eq:Kmax}) is true. For $\sigma>1$, condition (\ref{eq:V2}) does not yield any lower bound for $\mu$.
For $\sigma\in[1-d_p,1]$, condition (\ref{eq:V1}) is always true, since then $h(\sigma)\leq	\sigma\ln(\sigma)$. Finally, for $\sigma\in(0,1-d_p)$, condition (\ref{eq:U}) implies that the indicator functions in condition (\ref{eq:V1}) are equal to 1. Moreover, we have that
\begin{equation*}
	h(\sigma)\geq\sigma\ln(1-d_p)\,,
\end{equation*}
since the mapping $x\mapsto h(x)-x\ln(1-d_p)$ is decreasing/increasing on $(0,1)$ with minimum 0 at $x=1-d_p$. Combining these findings we have
\begin{equation*}
	u(\sigma)\,=\,\begin{cases}
		\qquad h(\sigma)\,,\quad&\sigma\in(0,1-d_p)\\[1ex]
		\sigma\ln(1-d_p)\,,\quad&\sigma\in[1-d_p,\infty)
	\end{cases}\,.
\end{equation*}

Now, let us derive the upper bound $o(\sigma)$ for $\mu$ in dependence of $\sigma$, such that the right-hand side of equivalence (\ref{eq:Kmax}) is valid. For $\sigma\leq1$, condition (\ref{eq:V1}) does not yield any upper bound for $\mu$. For $\sigma\in(1,1/(1-d_p)]$, it holds that $h(\sigma)\geq\ln(\sigma)$ such that condition (\ref{eq:V2}) is always true. Finally, for $\sigma>1/(1-d_p)$, condition (\ref{eq:U}) implies that the indicator functions in condition (\ref{eq:V2}) are equal to 1. Moreover, $h$ is increasing/decreasing on $(1,\infty)$ with maximum $-\ln(1-d_p)$ at $x=1/(1-d_p)$.  In summary, we have
\begin{equation*}
	o(\sigma)\,=\,\begin{cases}
		-\ln(1-d_p)\,,\quad&\sigma\in(0,1/(1-d_p)]\\[1ex]
		\qquad h(\sigma)\,,\quad &\sigma\in(1/(1-d_p),\infty)
	\end{cases}\,.
\end{equation*}

According to formula (\ref{eq:KSTpivot}) and as in the proof of Theorem \ref{thm:BandKST}, replacing $\bvt$ by $((\hat{\mu}-\mu)/\sigma,\hat{\sigma}/\sigma)$ then gives the representation for $C'_4$, from which the one for $C''_4$ is evident.

\bibliographystyle{tfnlm}  
\bibliography{BedburBib}

\end{document}